\documentclass[12pt]{article}
\usepackage{amsfonts,amsmath,amssymb,bbm}
\usepackage{setspace} 
\def\proof{\noindent{\bf Proof:}\hskip10pt}        
\def\QED{\hfill $\Box$ \vskip 2mm }

\font\tenmath=msbm10 scaled 1200
\font\sevenmath=msbm7 scaled 1200
\font\fivemath=msbm5 scaled 1200
\newfam\mathfam \textfont\mathfam=\tenmath
\scriptfont\mathfam=\sevenmath \scriptscriptfont\mathfam=\fivemath

\begin{document}
\def \\ { \cr }
\def\R{\mathbb{R}}
\def \1{1 \mkern -6mu 1} 
\def\N{\mathbb{N}}
\def\E{\mathbb{E}}
\def\P{\mathbb{P}}
\def\Z{\mathbb{Z}}
\def\Q{\mathbb{Q}}
\def\C{\mathbb{C}}
\def\D{\mathbb{D}}
\def\T{\mathbb{T}}
\def \e{{\rm e}}
\def \i{{\rm i}}
\def \f{{\mathcal F}}
\def \p{{\rm Part}}
\def \PN{{\rm Part}_{\N}}
\def \Pn{{\rm P}_{[n]}}
\def \g{{\mathcal G}}
\def \h{{\mathcal H}}
\def \d{{\tt d}}
\def \k{{\tt k}}
\def \r{{\mathcal R}}
\newtheorem{theorem}{Theorem}
\newtheorem{definition}{Definition}
\newtheorem{proposition}{Proposition}
\newtheorem{lemma}{Lemma}
\newtheorem{corollary}{Corollary}
\centerline{\LARGE \bf Asymptotic regimes for the occupancy}
\vskip 2mm
\centerline{\LARGE \bf  scheme of multiplicative cascades}
\vskip 1cm
\centerline{\Large \bf Jean Bertoin}
\vskip 1cm
\noindent
\centerline{\sl Laboratoire de Probabilit\'es, Universit\'e Pierre et Marie Curie}
\centerline{\sl  and DMA, Ecole Normale Sup\'erieure}
\centerline{\sl 175, rue du Chevaleret} 
\centerline{\sl F-75013 Paris, France}
\vskip 15mm

\noindent{\bf Summary. }{\small  In the classical occupancy scheme, one considers a fixed discrete probability measure  ${\bf p}=(p_i:  {i\in{\cal I}})$ and 
 throws balls independently at random in boxes labeled by ${\cal I}$, such that $p_i$ is the probability that a given ball falls into the box $i$. In this work, we are interested in asymptotic regimes of this scheme in the situation induced by a refining sequence $({\bf p}(k) : k\in\N)$ of random probability measures which arise from some multiplicative cascade.
Our motivation comes from the study of the asymptotic behavior of certain fragmentation chains.} 
 \vskip 5mm
\noindent
{\bf A.M.S. Classification.}  60 F 15, 60 J 80.
\vskip 3mm
\noindent{\bf e-mail.} {\tt jbe@ccr.jussieu.fr}

\begin{section}{Introduction}
 The  occupancy scheme is a simple urn model in probability theory that
 possesses a variety of applications to statistics, combinatorics, and computer science. These include, for instance, species sampling \cite{Bun, HeG}, analysis of algorithms \cite{Gardy}, learning theory \cite{Bou}, etc. The books by  Johnson and Kotz \cite{JK}  and by  Kolchin {\it et al.}   \cite{Kolchin} are standard references.
 
 This model is often depicted as  balls-in-bins. Typically, we denote by ${\rm Prob}_{\cal I}$
the space of probability measures on some countable set of indices ${\cal I}$, so  each 
${\bf p}\in{\rm Prob}_{\cal I}$ can be identified as a family  ${\bf p}=(p_i:  {i\in{\cal I}})$ of nonnegative real numbers with $\sum_{i\in{\cal I}} p_i = 1$. Given some ${\bf p}\in{\rm Prob}_{\cal I}$,  one throws balls successively and independently in a fixed series of boxes labeled by  indices $i$ in ${\cal I}$, and assumes that each ball has probability $p_i$ of falling into
the  box $i$. For every integers $j,n$ with $j\leq n$, we denote by ${N}_{n,j}^{\bf p}$ the number of boxes containing exactly $j$ balls when $n$ balls have been thrown, and by 
$${N}_{n}^{\bf p}:=\sum_{j=1}^{\infty}{N}_{n,j}^{\bf p}$$ the total number of occupied boxes.

We consider here a variant of this occupancy scheme which corresponds to
 a  nested family of boxes. This is conveniently described in terms of the genealogical structure of populations, so we start by recalling some notions in this area. We introduce the infinite {\it genealogical  tree}
$${\cal T}\,:=\,\bigcup_{k=0}^{\infty}\N^k\,,$$ with $\N:=\{1,2,\ldots\}$ and 
the convention
$\N^0:=\{\varnothing\}$. 
The elements of ${\cal T}$
are called {\it individuals}, and for every integer $k$, the $k$-th generation of 
${\cal T}$ is formed by the individuals in $\N^k$.
The boundary $\partial {\cal T}=\N^{\N}$ of ${\cal T}$
is the set of infinite sequences $\ell=(\ell_1, \ell_2, \ldots)$ of positive integers, which we call {\it leaves}. For each leaf $\ell=(\ell_1, \ell_2, \ldots)$ and each integer $k$, we write $\ell^{(k)}=(\ell_1,\ldots, \ell_k)$ for the ancestor of $\ell$ at generation $k$.
Conversely, for every individual  at generation $k$, say $i\in\N^k$, we denote by $\partial {\cal T}_i$
the subset of leaves whose ancestor at generation $k$ is $i$.
In particular, the root $\varnothing$ of the genealogical tree should be viewed as the progenitor of the entire population, and $\partial {\cal T}_{\varnothing}=\partial {\cal T}$.

Then consider some probability measure, say ${\bf P}$,  on $\partial {\cal T}$, and imagine that we sample a sequence $\lambda_{(1)},\lambda_{(2)}, \ldots$ of i.i.d. random leaves according to the law ${\bf P}$. For every fixed integers  $k,n\in\N$, we denote by
${N}^{(k)}_n$ the number of ancestors at generation $k$ of the first $n$ leaves :
\begin{eqnarray*}{N}^{(k)}_n
&:=&{\rm Card}\left\{\lambda_{(m)}^{(k)}: m\leq n\right\}\\
&=&{\rm Card}\left\{ i\in\N^k: 
{\rm Card}\left( \partial {\cal T}_i 
\cap \{\lambda_{(1)},\ldots,
\lambda_{(n)}\}\right)\geq 1\right\}\,.
\end{eqnarray*}
More precisely, we may also consider for every integer $1\leq j\leq n$
$${N}^{(k)}_{n,j}:
={\rm Card}\left\{ i\in\N^k: {\rm Card}\left( \partial {\cal T}_i 
\cap \{\lambda_{(1)},\ldots,
\lambda_{(n)}\}\right)= j\right\}\,,
$$
the number of individuals $i$ at generation $k$ such that the boundary $\partial {\cal T}_i$ of the subtree that stems from $i$ contains exactly $j$ leaves among $\{\lambda_{(1)},\ldots,
\lambda_{(n)}\}$.
The connexion with  the classical occupancy scheme may be better understood by viewing the random leaves as balls which are thrown on the boundary of the tree, and then imagining that each ball falls down following the branch from the leaf to the root $\varnothing$. Each individual $i$ can be thought of as a box, and if balls are thrown randomly according to the probability measure ${\bf P}$ on the boundary of the tree, then the probability that some given ball passes through the box $i$ at generation $k$ is  
$$p_i(k):={\bf P}(\partial {\cal T}_i)\,.$$
Clearly, ${\bf p}(k):=(p_i(k):  i\in \N^k)$  defines a probability measure on $\N^k$ for each generation $k$, and the sequence of discrete probability measures $({\bf p}(k) : k\in\N)$ determines ${\bf P}$.

Recently, Gnedin {\it et al.} \cite{GPY} have considered asymptotic laws for a randomized version of the classical occupation scheme, where the discrete probability ${\bf p}\in{\rm Prob}_{\cal I}$ is random (more precisely, ${\bf p}$ is obtained from the atoms of some Poisson point measure on $]0,\infty[$). 
In the present work, we will be interested in a situation when the probability measure ${\bf P}$ on $\partial {\cal T}$ (and therefore also each probability measure ${\bf p}(k)$ on $\N^k$) is random. So henceforth, {\it conditionally on ${\bf P}$},  each leaf $\lambda_{(1)}, \ldots$ is picked randomly according to ${\bf P}$ and independently of the others. 
This is equivalent to assuming that the sequence of random leaves is exchangeable
with de Finetti measure ${\bf P}$. More precisely, we shall assume that ${\bf P}$ is given by some {\it multiplicative cascade}; see Liu \cite{Liu} and the references therein.
This means that we consider first some random probability measure $\pmb{\varrho}=(\varrho_1, \ldots)$ 
in ${\rm Prob}_{\N}$ and assign to each individual $i=(i_1, \ldots, i_k)$ of the genealogical tree an independent copy $\pmb{\varrho}(i)$ of $\pmb{\varrho}$. Roughly speaking, 
$\pmb{\varrho}(i)$ describes how the mass $p_i(k)={\bf P}(\partial {\cal T}_i)$ is splitted to
the subsets of leaves $\partial {\cal T}_{ij}$ for $j\in\N$, where
$ij=(i_1, \ldots, i_k,j)$ denotes the $j$-th child of the individual $i$ at generation $k+1$. Specifically, $\varrho_j(i)$ is the portion of the mass of $i$ inherited by the child $ij$, i.e. 
$$\frac{p_{ij}(k+1)}{p_i(k)}=\varrho_j(i)\,,$$
so that, by iteration,
\begin{equation}\label{eq3}
p_i(k)=\varrho_{i_1}(\varnothing)\times \varrho_{i_2}(i^{(1)})\times \cdots \times \varrho_{i_k}(i^{(k-1)})
\end{equation}
where $i^{(k')}=(i_1, \ldots, i_{k'})$ denotes the ancestor of $i$ at generation $k'\leq k$. 
Clearly,  for each integer $k$, ${\bf p}(k)=(p_i(k):  i\in \N^k)$ now defines a random probability measure on $\N^k$, and we can identify the conditional laws
\begin{equation}\label{eqidl}
{\cal L}\left({N}^{(k)}_n \mid {\bf p}(k)\right) \ = \
{\cal L}\left({N}^{{\bf p}(k)}_n\right)
\ \hbox{ and }\  
{\cal L}\left({N}^{(k)}_{n,j} \mid {\bf p}(k)\right) 
\ =\ 
{\cal L}\left({N}^{{\bf p}(k)}_{n,j}\right)\,.
\end{equation}

Our main purpose is to determine the asymptotic regimes of
the numbers of occupied boxes  ${N}^{(k)}_{n,j}$ and ${N}^{(k)}_{n}$ when both $n$ and $k$ tend to infinity.
It is easily seen from routine estimates that non-degenerate limits should occur when $k\approx \ln n$. Since both $k$ and $n$ are integers, a natural regime thus could be
$k=\lfloor a\ln n\rfloor$ for some real number $a>0$, where the notation $\lfloor\cdot\rfloor$ refers to the integer part. It turns out that this is actually too crude, in the sense that for
$k=\lfloor a\ln n\rfloor$, the asymptotic behavior of ${N}^{(k)}_{n,j}$
does not only depend on $a$, but also on the oscillations
of the fractional part $\{a\ln n\}$.  
Indeed, we shall establish a law of large numbers  for ${N}^{(k)}_{n,j}$
and a central limit theorem for ${N}^{(k)}_{n}$
 when  $k,n\to \infty$ in such a way that $k=a \ln n+ b + o(1)$ for fixed real numbers $a$ and $b$ in certain intervals.

Our approach essentially combines 
uniform probability estimates for the classical occupancy scheme
and information about asymptotic behaviors in multiplicative cascades which can be gleaned from the literature and  will be reviewed in Section 2.
In particular,  the analysis of multiplicative cascades relies crucially 
on the natural connexion with a class of branching random walks, and more precisely, on their large deviations behaviors whose descriptions are  due to Biggins \cite{Big2}.  
The main results about the asymptotic regimes in our model will presented and proved
in Section 3. They include a law of large numbers and a central limit theorem mentioned above; we will also study asymptotics of the shattering generation, i.e. the lowest generation at which no box contains more than a fixed number of balls.
Finally, we shall conclude this work by discussing some interpretations of the present results in the framework of homogeneous fragmentation processes, which provided the initial motivation for this work.

\end{section}

\begin{section}{Preliminaries}

\subsection{Some limit theorems for the occupancy scheme}
In this section, we lift from the literature on urn models  a law of large numbers and a central limit theorem for the number of occupied boxes
that will be useful in our study.

Given an arbitrary discrete probability measure  ${\bf p}=(p_i:  {i\in{\cal I}})\in{\rm Prob}_{\cal I}$,  we define first for every $j\leq n$ 
the number
$\bar{N}^{\bf p}_{n,j}$  of boxes occupied by more than $j$ balls when $n$ balls have been thrown, viz. 
$$\bar{N}^{\bf p}_{n,j}:=\sum_{{\ell}=j+1}^n{N}^{\bf p}_{n,{\ell}}\,.$$
In particular, for $j=0$,  $\bar{N}^{\bf p}_{n,0}= {N}^{\bf p}_{n}$.
Introduce also
for every real number $x\geq 0$
\begin{eqnarray*}
\bar\mu^{\bf p}_{j}(x)&:=&\sum_{{\ell}=j+1}^{\infty}\frac{x^{\ell}}{{\ell}!}\sum_{i\in{\cal I}} p_i^{\ell}\e^{-p_i x}\\
&=&\sum_{i\in{\cal I}} \left(1-\e^{-p_i x}\left(1+\cdots
+\frac{(p_ix)^{j}}{j!}\right)\right)\,.
\end{eqnarray*}
We may now state the following  law of large numbers which has its roots in 
Bahadur \cite{Bahadur}.

 \begin{lemma}\label{lln} Let ${\bf p}(1), {\bf p}(2), \ldots$ be a sequence  of discrete probability measures,  $(n_k, k\in \N)$ a sequence of positive integers with
 $\lim_{k\to\infty}n_k=\infty$,
  and $j\in\Z_+$. Suppose that
\begin{equation}\label{BC}
\sum_{k\in\N} \frac{1}{\bar\mu^{{\bf p}(k)}_{j}(n_k)}<\infty
\end{equation}
 and
 $$\lim_{\alpha\to 1}\lim_{k\to\infty} \frac{\bar\mu^{{\bf p}(k)}_{j}(\alpha n_k)}
 {\bar\mu^{{\bf p}(k)}_{j}(n_k)}=1\,.$$
 Then
 $$\lim_{k\to\infty}\frac{\bar{N}_{n_k,j}^{{\bf p}(k)}}{ \bar\mu^{{\bf p}(k)}_{j }(n_k)}=1\qquad \hbox{a.s.} $$
   \end{lemma}
Although this result should belong to the folklore of limit theorems for urn models,
we have not been able to find a precise reference where it is stated in this form,
and thus we shall provide a proof. The argument relies on  Poissonization, which is an important technique in this area; see, for instance, the surveys by Gnedin {\it et al.} \cite{GHP} or Holst \cite{Holst}.

\proof We work first with a fixed
 probability measure  ${\bf p}=(p_i:  {i\in{\cal I}})$, but we replace the 
deterministic number of balls $n$ by ${\tt n}_x$, where 
${\tt n}=({\tt n}_x, x\geq 0)$ is an independent
standard Poisson process.  The key effect of Poissonization is that now, for each $i\in{\cal I}$, the number of balls in the box $i$ has the Poisson distribution with parameter $p_i x$,  and that to different boxes correspond independent Poisson variables. As a consequence,  the variable $\beta_i$ which takes the value $1$ 
if more than $j$ balls occupy the box $i$ and $0$ otherwise, has the Bernoulli distribution
with parameter
$$\sum_{{\ell}=j+1}^{\infty}\e^{-p_i x}\frac{(p_i x)^{\ell}}{{\ell}!}\,,$$
and when $i$ varies in ${\cal I}$,  these Bernoulli variables are independent.
Changing the typography, we write 
$$\bar{\tt N}^{\bf p}_{x,j}:=\sum_{i\in{\cal I}}\beta_i$$
 for the number of boxes occupied by more than $j$ balls when ${\tt n}_x$ balls have been thrown. By elementary properties of sums of independent Bernoulli variables, we see that
\begin{equation}\label{eq45}
\E\left(\bar{\tt N}^{\bf p}_{x,j}\right)=\bar\mu^{\bf p}_{j}(x)
\end{equation}
and 
$${\rm Var}\left(\bar{\tt N}^{\bf p}_{x,j}\right)\leq \bar\mu^{\bf p}_{j}(x)\,.$$
Thus Chebyshev's inequality ensures that 
$$\P\left( \left|\frac{ \bar{\tt N}^{\bf p}_{x,j}}{\bar\mu^{\bf p}_{j}(x)} -1\right|\geq \varepsilon  \right)
\leq \frac{\varepsilon^{-2}}{\bar\mu^{\bf p}_{j}(x)}$$
for every $\varepsilon>0$.

Next, we replace the fixed probability measure  ${\bf p}$ by ${\bf p}(k)$ and take $x=\alpha n_k$ for some real number $\alpha$ close to $1$. The bound above combined our assumptions enables us to apply the Borel-Cantelli lemma, and we get that
\begin{equation}\label{eq6}
\lim_{k\to\infty}\frac{\bar{\tt N}_{\alpha n_k,j}^{{\bf p}(k)}}{ \bar\mu^{{\bf p}(k)}_{j }(n_k)}={\ell}(\alpha)
\qquad \hbox{a.s.}\,,
\end{equation}
with
$${\ell}(\alpha):=\lim_{k\to\infty} \frac{\bar\mu^{{\bf p}(k)}_{j}(\alpha n_k)}
 {\bar\mu^{{\bf p}(k)}_{j}(n_k)}\,.$$

Recall that the  number ${\tt n}_x={\tt n}_{\alpha n_k}$ of balls which are thrown has the Poisson
distribution with parameter $\alpha n_k$. On the event $\left\{{\tt n}_{\alpha n_k}\geq n_k\right\}$, there is the 
bound 
$$\bar{\tt N}_{\alpha n_k,j}^{{\bf p}(k)}\geq \bar{ N}_{n_k,j}^{{\bf p}(k)}\,,$$
whereas on the complementary event we have
$$\bar{\tt N}_{\alpha n_k,j}^{{\bf p}(k)}\leq \bar{ N}_{n_k,j}^{{\bf p}(k)}\,.$$
Recall that $n_k\to\infty$.
Plainly, if $\alpha>1$, then 
$$\lim_{k_0\to \infty}\P({\tt n}_{\alpha n_k}\geq n_k\hbox{ for all }k\geq k_0) = 1$$
whereas  if $\alpha<1$, then 
$$\lim_{k_0\to \infty}\P({\tt n}_{\alpha n_k}< n_k\hbox{ for all }k\geq k_0) = 1.$$
This completes the proof, by using \eqref{eq6} and the assumption that $\lim_{\alpha\to1}{\ell}(\alpha)=1$. \QED

\noindent {\bf Remark.} If we replace the requirement \eqref{BC} in Lemma \ref{lln}
by the weaker $\lim_{k\to\infty} \bar\mu^{{\bf p}(k)}_{j}(n_k)=\infty$,
the same calculations yield the weak law of large numbers :
 $$\lim_{k\to\infty}\frac{\bar{N}_{n_k,j}^{{\bf p}(k)}}{ \bar\mu^{{\bf p}(k)}_{j }(n_k)}=1\qquad \hbox{in probability.} $$

\vskip 2mm

Next, we turn our attention to fluctuations for the number of occupied boxes.
Following Hwang and Janson \cite{HJ}, we introduce for every  fixed probability measure  ${\bf p}=(p_i:  {i\in{\cal I}}) \in {\rm Prob}_{\cal I}$ and every
$x\geq 0$
\begin{equation}\label{meanHJ}
\mu_{\bf p}(x):=\bar\mu^{\bf p}_{0}(x)=\sum_{i\in{\cal I}}(1-\e^{-p_i x})
\end{equation}
and 
\begin{equation}\label{varHJ}
\sigma^2_{\bf p}(x)
:=\sum_{i\in{\cal I}}\e^{-p_i x}(1-\e^{-p_i x})-
 x^{-1}\left(\sum_{i\in{\cal I}}xp_i \e^{-p_i x}\right)^2\,.
 \end{equation}

These quantities provide uniform estimates for the mean and the variance of
 ${N}_{n}^{\bf p}$; specifically it is known  from Theorem 2.3 in \cite{HJ} that 
 \begin{equation}\label{eq4}
 |\E({N}_{n}^{\bf p}) - \mu_{\bf p}(n)| \leq c
 \end{equation}
 and 
   \begin{equation}\label{eq4'}
   |{\rm Var}({N}_{n}^{\bf p}) - \sigma^2_{\bf p}(n)|\leq c\,, 
  \end{equation}
where $c$ denotes some numerical constant (which depends neither of $n$ nor of ${\bf p}$). This makes the following central limit theorem quite intuitive
(see Corollary 2.5 in \cite{HJ}, and also Dutko \cite{Dutko} and Karlin \cite{Karlin} for earlier versions).
 
 \begin{lemma}\label{LHJ} Let ${\bf p}(1), {\bf p}(2), \ldots$ be a sequence of discrete probability measures and $(k_n, n\in \N)$ a sequence of positive integers such that
 $$\lim_{n\to\infty}\sigma^2_{{\bf p}(k_n)}(n)=\infty\,.$$
 Then the number of occupied boxes is asymptotically normally distributed when $n$ goes to infinity, in the sense that
 $$ \frac{{N}_{n}^{{\bf p}(k_n)} - \mu_{{\bf p}(k_n)}(n) }{\sigma_{{\bf p}(k_n)}(n)}$$
 converges in distribution to a standard normal variable as $n\to \infty$.
\end{lemma}
We do not know whether a similar central limit theorem holds for the number
${N}^{\bf p}_{n,j}$  of boxes occupied by exactly $j$ balls for $j\geq 1$.

\subsection{Large deviations behaviors of multiplicative cascades}

Recall that $\pmb{\varrho}$ is a random probability measure on $\N$. We
denote its law by $\nu$, so $\nu$ is a probability measure on ${\rm Prob}_{\N}$
that will be referred to as the {\it splitting law}.
We shall always assume that this splitting law is not geometric
\footnote{Working with a geometric splitting law would induce a phenomenon of periodicity which we shall not discuss here for simplicity. However results similar to those proven in this work can be established by the same techniques for geometric splitting  laws.}, in the sense that
there is no real number $r>0$ such that with probability one, all the atoms
of $\pmb{\varrho}$ belong to $\{r^n, n\in\Z_+\}$.
In particular, note that the degenerate case when $\pmb{\varrho}$ is a Dirac point mass a.s. is henceforth excluded.

As it was explained in the Introduction, we consider a family $(\pmb{\varrho}(i), i\in{\cal T})$ of independent copies of $\pmb{\varrho}$ labeled by the individuals of the genealogical tree ${\cal T}$. 
The multiplicative cascade construction \eqref{eq3} defines a random probability measure ${\bf p}(k)$ on $\N^k$ for every generation $k\in\N$.
Our aim is to apply general asymptotic results for occupancy schemes such as Lemmas  \ref{lln} and \ref{LHJ}, and in this direction, we shall use fundamental large deviations behaviors  for branching random walks that  Biggins \cite{Big2} established.

Taking logarithm of masses, we may encode the random probability measure
${\bf p}(k)=(p_i(k): i\in\N^k)$ at generation $k$ by
 the random point measure on $\R_+$
$$Z^{(k)}({\rm d}y) := \sum \delta_{-\ln p_i(k)}({\rm d}y)\,,$$
where $\delta_z$ stands for the Dirac point mass at $z$ and the sum in the right-hand side is taken over the individuals $i$ at the $k$-th generation which have a positive mass.
It then follows immediately from the structure of multiplicative cascade \eqref{eq3} that
$(Z^{(k)}, k\in\Z_+)$ is a branching random walk, in the sense that
for every integers $k,k'\geq 0$, $Z^{(k+k')}$ is obtained from
$Z^{(k)}$ by replacing each atom $z$ of $Z^{(k)}$ by
a family $\{z+y, y\in{\cal Y}\}$, where ${\cal Y}$
is distributed as the family of the atoms of $Z^{(k')}$ and distinct atoms $z$ of $Z^{(k)}$ correspond to independent copies of ${\cal Y}$.

We now introduce analytic quantities defined in terms of the splitting law $\nu$
which will have an important role in the present study.
First, we define the Laplace transform of the intensity measure of $Z^{(1)}$ by
$$
\hbox{\sc l}(\theta):=\E(\langle Z^{(1)}, \e^{-\theta \cdot}\rangle)
$$
 for $\theta >0$; note that there are also the alternative expressions  
\begin{equation}\label{defm}
\hbox{\sc l}(\theta)=\E\left(\sum_{j\in\N} \varrho_j^{\theta}\right)
\,=\,\int_{{\rm Prob}_{\N}}
\left( \sum_{i\in\N}p_i^{\theta}\right)\nu({\rm d}{\bf p})\,.
\end{equation}
The function $\hbox{\sc l}: ]0,\infty[\to ]0,\infty]$ is convex decreasing with $\hbox{\sc l}(1)=1$; we define
\begin{equation}\label{defstar}
\theta_*\,:=\,\inf\left\{\theta >0: \hbox{\sc l}(\theta)<\infty\right\}\,,
\end{equation} 
so that $\hbox{\sc l}(\theta)<\infty$ when $\theta>\theta_*$.
One readily sees from H\"older's inequality that  $\ln \hbox{\sc l}$ is a convex function,
and then that
\begin{equation}\label{defvarphi}
\varphi (\theta):=\ln \hbox{\sc l}(\theta)-\theta \frac{\hbox{\sc l}'(\theta)}{\hbox{\sc l}(\theta)}
\end{equation}
is a function which decreases on $]\theta_*,\infty[$.

As $\hbox{\sc l}(1)=1$ and $\hbox{\sc l}$ decreases, we have $\varphi(1)=-\hbox{\sc l}'(1)>0$, and thus the set of $\theta\in ]\theta_*,\infty[$ such that $\varphi(\theta)>0$ is a non-empty open interval $]\theta_*,\theta^*[$, where  $ \theta^*>1$ is defined by
\begin{equation}\label{defetoile}
\theta^*:=\sup\{\theta >\theta_* : \varphi(\theta)>0\}\,.
\end{equation}

\noindent{\bf Remark :} The critical parameter $ \theta^*$ may be finite or infinite, and is finite whenever  $$\Vert \max_{i\in \N} \varrho_i\Vert_{\infty}=1\,,$$
where $\pmb{\varrho}=(\varrho_i, i\in\N)$ denotes  a random probability measure on $\N$ with law $\nu$. Indeed, it is easily seen that
$$\lim_{\theta\to\infty} \hbox{\sc l}(\theta)^{1/\theta} = \Vert \max_{i\in \N} \varrho_i\Vert_{\infty}\,,$$
and when the right-hand side equals $1$, the function $g: \theta \to 
 -\frac{\ln  \hbox{\sc l}(\theta)}{\theta}$ has thus limit $0$ at infinity. Since $g$ is non-negative on $[1,\infty[$ and $g(1)=0$, $g$ reaches its overall maximun at some location at, say,  $\theta_{\rm max}\in ]1,\infty[$. As $g'(\theta)=\theta^{-2}\varphi(\theta)$, we conclude that $\theta_{\rm max}=\theta^*<\infty$.

Following Biggins \cite{Big1}, we are now able to introduce
for every $\theta>\theta_*$
$$W^{(k)}(\theta):=\hbox{\sc l}(\theta)^{-k}\langle Z^{(k)},\e^{-\theta \cdot}\rangle
= \hbox{\sc l}(\theta)^{-k}\sum_{i\in\N^k}p_i(k)^{\theta}\,,\qquad k\geq 0\,,$$
which form 
 a remarkable family of martingales :

\begin{lemma}\label{L1} For every $\theta\in]\theta_*,\theta^*[$,
 the martingale $(W^{(k)}(\theta), k\in \Z_+)$ is bounded in $L^{\gamma}(\P)$ for some $\gamma>1$. Its terminal value  
 $$W(\theta):=\lim_{k\to\infty}W^{(k)}(\theta)$$
 is (strictly) positive a.s.
 \end{lemma}
 
\proof Jensen's inequality implies that for every probability measure 
${\bf p}\in{\rm Prob}_{\N}$ and every $\gamma>1$, there is the upper-bound
$$\left( \sum_{i\in\N}p_i^{\theta}\right)^{\gamma}
\leq \sum_{i\in\N}p_i^{\gamma(\theta-1)+1}.$$
For any $\theta>\theta_*$, we may chose $\gamma>1$
sufficiently small  
such that $\gamma(\theta-1)+1>\theta_*$, and we deduce that 
 $\E(W^{(1)}(\theta)^{\gamma})<\infty$.

We then 
observe that the function $f:\theta\to \theta^{-1}\ln \hbox{\sc l}(\theta)$ has derivative $f'(\theta)=-\theta^{-2} \varphi(\theta)$. 
Thus this derivative is negative when $\theta\in]\theta_*,\theta^*[$, which means that 
$f$ decreases in some neighborhood of $\theta$.
We may thus find $\gamma >1$ sufficiently small  such that
$$\frac{\ln \hbox{\sc l}(\gamma \theta)}{\gamma \theta} < \frac{\ln \hbox{\sc l}(\theta)}{\theta}\,,$$
and hence $\hbox{\sc l}(\gamma \theta)<\hbox{\sc l}(\theta)^{\gamma}$.
We can now apply Theorem 1 in Biggins \cite{Big2}, which completes the proof of the first part of our claim. Finally, the assertion that the terminal value $W(\theta)>0$ a.s. derives easily from the
fact the probability that branching random walk $Z^{(k)}$ is extinguished 
at generation $k$ equals $0$ for every $k\in\N$.  \QED

 In order to state the key technical result  for this present study, 
we define for every $\theta>\theta_*$ the tilted random point measure
$$Z^{(k)}_{\theta}({\rm d}y):=\frac{\e^{-\theta y}}{\hbox{\sc l}(\theta)^k}Z^{(k)}({\rm d} y)\,,\qquad y\geq 0\,.$$
We also  introduce the mean and the variance of the intensity measure of
$Z^{(1)}_{\theta}$ :
\begin{equation}\label{deftilde}
\hbox{\sc m}(\theta):=-\frac{\hbox{\sc l}'(\theta)}{\hbox{\sc l}(\theta)}\ \hbox{ and }\ 
\hbox{\sc v}(\theta)= \frac{\hbox{\sc l}''(\theta)}{\hbox{\sc l}(\theta)}-\left(\frac{\hbox{\sc l}'(\theta)}{\hbox{\sc l}(\theta)}\right)^2\,.
\end{equation}
Both $\hbox{\sc m}(\theta)$ and $\hbox{\sc v}(\theta)$ are positive quantities,
and write $g_{\theta}$ for the (centered) Gaussian density with variance $\hbox{\sc v}(\theta)$, i.e.
$$g_{\theta}(x)=\frac{1}{\sqrt{2\pi \hbox{\sc v}(\theta)}}\exp\left(-\frac{x^2}{2{\hbox{\sc v}(\theta)}}\right)\,.$$

The next statement is a version of Theorem 4 in Biggins \cite{Big2}
specialized to our framework.

\begin{lemma}\label{L2} 
 The following assertion holds with probability one:
$$\lim_{k\to\infty}
\left|\sqrt k Z^{(k)}_{\theta}([x+ k \hbox{\sc m}(\theta) -h, x+ k \hbox{\sc m}(\theta) +h[)-2hW(\theta) g_{\theta}(x/\sqrt k)\right| = 0\,,
$$
where the limit is uniform for $x\in\R$, $h\leq 1$ and $\theta$ in a compact subset of 
$]\theta_*,\theta^*[$.
\end{lemma}

Next, observe that if $f:\R\to\R_+$ is, say, a continuous function, then
\begin{eqnarray*}
\sum_{i\in\N^k}Êf(k \hbox{\sc m}(\theta) +\ln p_i(k)) &=& \int_{\R_+}Êf(k \hbox{\sc m}(\theta) -y) Z^{(k)}(  {\rm d}y)\\
&=& 
\left(\hbox{\sc l}(\theta)\e^{\theta\hbox{\sc m}(\theta)}\right)^k
\int_{\R_+}Êf(k \hbox{\sc m}(\theta) -y) \e^{\theta(y-k \hbox{\sc m}(\theta))} Z^{(k)}_{\theta}(  {\rm d}y).
\end{eqnarray*}
Recall also that the rate function $\varphi$ is defined by \eqref{defvarphi}. We finally state the following limit theorem  which will be useful to estimate the conditional mean number of occupied boxes given the multiplicative cascade.

\begin{corollary}\label{C1} {\rm (large deviations behavior)}
 Pick $\theta\in ]\theta_*,\theta^*[$ and let
$f:\R\to\R_+$ be a continuous function.
Assume that there exist $\alpha >0$ and $\beta>\theta$ such that
$$\lim_{y\to +\infty} y^{\alpha} f(y) =0 \quad \hbox{and} \quad
 \lim_{y\to -\infty} \e^{-\beta y} f(y) = 0\,,$$
so in particular $f\in L^1( \e^{-\theta y}{\rm d} y)$.
Let also $(c_k: k\in\N)$ denote a sequence of real numbers which converges to some $c\in\R$.
Then with probability one, we have
$$\lim_{k\to\infty}\sqrt{k }\e^{-\varphi(\theta)k}\sum_{i\in\N^k}Êf(k \hbox{\sc m}(\theta) +\ln p_i(k) + c_k)
= \frac{\e^{\theta c}}{\sqrt{2\pi \hbox{\sc v}(\theta)}} \left(\int_{\R} f(y)\e^{-\theta y}
{\rm d}y\right) W(\theta)\,.$$

\end{corollary}
Corollary \ref{C1} follows from readily from Lemma \ref{L2} when $f$ has bounded support. However, the derivation in the case when the function $f$ has unbounded 
support is rather technical, even though our assumptions have been tailored for the purpose of the present work.   We postpone the proof to the Appendix.

\end{section}

\begin{section}{Asymptotic regimes}

\subsection{Main results}

We now have all the technical ingredients for our study, we just need to introduce a few more notation. 
We shall consider the regime for pairs of integers $(k,n)$ such that  
\begin{equation}\label{reg} k,n\to\infty \quad \hbox{and}\quad
k-a\ln n \to b\,,
\end{equation}
where $a>0$ and $b\in\R$ are fixed. When $F(k,n)$ is some function depending on
$k$ and $n$, we shall write
$$\lim_{a,b}F(k,n)$$
for the limit of $F(k,n)$ when $(k,n)$ follows the regime \eqref{reg}, of course provided that such a limit exists.

Recall  that our basic datum is the splitting law $\nu$ on ${\rm Prob}_{\N}$, and that its Laplace transform $\hbox{\sc l}(\theta)$
is given by \eqref{defm}. Further important notions include the critical parameters 
$\theta_*,\theta^*$, the rate function $\varphi$, 
 and the mean $\hbox{\sc m}$ and variance $\hbox{\sc v}$ functions,
 which have been  defined in \eqref{defstar}, \eqref{defetoile},  \eqref{defvarphi} and \eqref{deftilde}, respectively.
The mean function $\hbox{\sc m}$ decreases continuously on $]\theta_*,\theta^*[$ and takes positive values. We denote the inverse bijection by
$$\hbox{\sc m}^{-1}: ]\hbox{\sc m}_*, \hbox{\sc m}^*[\to ]\theta_*,\theta^*[\,,$$
where
$$\hbox{\sc m}_* = \lim_{\theta\to \theta^*-}\hbox{\sc m}(\theta)
\quad \hbox{and}\quad 
\hbox{\sc m}^* = \lim_{\theta\to \theta_*+}\hbox{\sc m}(\theta)\,.$$
One always has $]\hbox{\sc m}_*, \hbox{\sc m}^*[\subseteq ]0,\infty[$, and the inclusion can be strict.

\noindent{\bf Example :}  These quantities are especially simple in the case when $\nu$ is the Poisson-Dirichlet distribution ${\rm PD}(1)$. Indeed, one easily  gets  $\hbox{\sc l}(\theta)=
1/ \theta$ for $\theta >0$,   and then $\varphi(\theta)
=-\ln \theta +1$. One thus sees that  $\theta_* = 0$ and $\theta^*= \e$. Finally 
 $\hbox{\sc m}(\theta)=1/\theta$ and $\hbox{\sc v}(\theta)=1/\theta^{2}$, so $\hbox{\sc m}_*=1/\e$,  $\hbox{\sc m}^*=\infty$ and  $\hbox{\sc m}^{-1}(a)=1/a$.

\vskip 2mm

Finally recall from Lemma \ref{L1} that for every $\theta\in]\theta_{*},\theta^{*}[$, $W(\theta)$ is a positive random variable which arises as the limit of a remarkable martingale.
We are now able to specify regimes for the (strong) law of large numbers for the number of boxes occupied by exactly $j$ balls in an occupation scheme driven by a multiplicative cascade.

\begin{theorem} \label{T1}
Pick $a\in]1/\hbox{\sc m}^*, 1/\hbox{\sc m}_*[$ and $b\in\R$, and set $\theta=\hbox{\sc m}^{-1}(1/a)$. 
Then for every integer $j>\theta$, the following limits
$$\lim_{a,b} \sqrt{k }\e^{-\varphi(\theta)k} \, 
\bar N^{(k)}_{n,j-1}
=\left(\sum_{\ell=j}^{\infty}\frac{\Gamma(\ell-\theta)}{\ell !}\right)
\frac{\e^{-\theta b/a}}{\sqrt{2\pi \hbox{\sc v}(\theta)}}\,W(\theta)\,,$$
and 
$$\lim_{a,b} \sqrt{k} \e^{-\varphi(\theta)k}\, 
 N^{(k)}_{n,j}
=\frac{\Gamma(j-\theta)\e^{-\theta b/a}}{j !\sqrt{2\pi \hbox{\sc v}(\theta)}}\,W(\theta)\,,$$
hold with probability one.
\end{theorem} 

\noindent{\bf Remark : } It may be interesting to recall that $ \frac{1}{\sqrt k}\e^{\varphi(\theta)k}$ is also the order of magnitude of the numbers of boxes of size approximately 
$\e^{-k/a}\approx 1/n$ at generation $k$; see Corollary 3 in \cite{BeRou}. We further stress that  $\varphi(\theta)\leq 1/a$ and that this inequality is strict except when $\theta=1$
(this can be checked directly from the observation that   $\ln  \hbox{\sc l}$ is a strictly convex function).

\proof 
We work conditionally on the random probabilities 
${\bf p}(k)$ using \eqref{eqidl}. We aim at applying Lemma \ref{lln}, and in this direction we fix some real number $\alpha$ close to $1$ and observe that
\begin{eqnarray*}
\bar\mu^{{\bf p}(k)}_{j-1}(\alpha n)&=&\sum_{i\in\N^{k}} \left(1-\e^{-p_i(k) \alpha n}\left(1+\cdots
+\frac{(p_i(k) \alpha n)^{j-1}}{(j-1)!}\right)\right)\\
&=&
\sum_{i\in\N^{k}}Êf(k/a +\ln p_i(k) + c_{k} )
\end{eqnarray*}
with
$$f(x)=  \left(1-\left(1+\cdots
+\frac{\e^{(j-1)x}}{(j-1)!}\right)\exp\{-\e^x\}\right)$$
and 
$$c_{k}:= \ln n + \ln \alpha - k  /a = -b/a+ \ln \alpha + o(1)\,.$$

We can now check that the assumptions of Corollary \ref{C1} hold with 
 $c=-b/a+ \ln \alpha$. A straightforward calculation shows that
 $$\int_{\R} f(y)\e^{-\theta y}
{\rm d}y=  \sum_{\ell=j}^{\infty}\frac{\Gamma(\ell-\theta)}{\ell !}\,,$$
and then, invoking Corollary \ref{C1} and  recalling that $\theta=\hbox{\sc m}^{-1}(1/a)$, we deduce 
\begin{equation}\label{estmean}
\lim_{a,b} \sqrt{k }\e^{-\varphi(\theta)k} \, 
\bar\mu^{{\bf p}(k)}_{j-1}(\alpha n)
=\left(\sum_{\ell=j}^{\infty}\frac{\Gamma(\ell-\theta)}{\ell !}\right)
\frac{ \e^{\theta (-b/a+ \ln \alpha)}}{\sqrt{2\pi \hbox{\sc v}(\theta)}}\,W(\theta)\,.
\end{equation}

Next, we fix $\eta\in[-1,1]$ and define for every integer $k$
$$n_{k,\eta}:=\left \lfloor\exp\left((k-b-\eta)/a\right)\right\rfloor\,.$$
Replacing $b$ by $b+\eta$ and $n$ by $n_{k,\eta}$  in  \eqref{estmean},
we get
$$\lim_{k\to\infty}\ \sqrt{k } \e^{-\varphi(\theta)k}\, 
\bar\mu^{{\bf p}(k)}_{j-1}(\alpha n_{k,\eta})
=\left(\sum_{\ell=j}^{\infty}\frac{\Gamma(\ell-\theta)}{\ell !}\right)
\frac{ \e^{\theta ( \ln \alpha -(b+\eta)/a)}}{\sqrt{2\pi \hbox{\sc v}(\theta)}}\,W(\theta)\,.$$
Note that the right-hand side  depends continuously on the variable $\alpha$
and that the requirement \eqref{BC} in  Lemma \ref{lln} is fulfilled, as $\varphi(\theta)>0$.
An application of this lemma gives 
$$\lim_{k\to\infty}\ \sqrt{k }e^{-\varphi(\theta)k} \, 
\bar N^{(k)}_{n_{k,\eta},j-1}
=\left(\sum_{\ell=j}^{\infty}\frac{\Gamma(\ell-\theta)}{\ell !}\right)
\frac{\e^{-\theta (b+\eta)/a}}{\sqrt{2\pi \hbox{\sc v}(\theta)}}\,W(\theta)\,,$$
with probability one. An argument of monotonicity, namely
$$n_{k,\eta}Ê\leq n \leq n_{k,\eta'}Ê\ \Longrightarrow \ 
\bar N^{(k)}_{n_{k,\eta},j-1}\leq \bar N^{(k)}_{n,j-1}\leq \bar N^{(k)}_{n_{k,\eta'},j-1}\,,$$
 completes the proof of our first claim. The second follows immediately from the first, since  $N^{(k)}_{n,j}=\bar N^{(k)}_{n,j-1}-\bar N^{(k)}_{n,j}$. 
\QED

We next turn our attention to finer asymptotic results  for the total number of occupied boxes 
 $\bar N^{(k)}_{n,0}= N^{(k)}_{n}$. If $\theta<1$, that is if 
$a< 1/\hbox{\sc m}(1)$, then we can take $j=1$ in Theorem \ref{T1},
and we get from the easy identity 
$$ \sum_{\ell=1}^{\infty}\frac{\Gamma(\ell-\theta)}{\ell !}=
\int_0^{\infty}(1-\e^{-x})x^{-\theta-1}{\rm d}x=\frac{\Gamma(1-\theta)}{\theta} 
$$
that
\begin{equation}\label{eq67}
\lim_{a,b} \sqrt{k}\e^{-\varphi(\theta)k}\, N^{(k)}_n
=\frac{\Gamma(1-\theta)\e^{-\theta b/a}}{\theta \sqrt{2\pi \hbox{\sc v}(\theta)}}\, W(\theta)
\qquad \hbox{a.s.}\end{equation}
Recall also from \eqref{eq4} that the conditional expectation of number of occupied boxes given the de Finetti measure, $\E(N^{(k)}_n\mid {\bf p}(k))$,
can be approximated by $\mu_{{\bf p}(k)}(n)=\bar\mu^{{\bf p}(k)}_{0}(n)$,
and that \eqref{estmean} provides an estimation of the latter. 
This gives
\begin{equation}\label{eq68}
\lim_{a,b} \sqrt{k}\e^{-\varphi(\theta)k}\, \mu_{{\bf p}(k)}(n)=\frac{\e^{-\theta b/a}\Gamma(1-\theta)}{\theta \sqrt{2\pi \hbox{\sc v}(\theta)}} W(\theta)\,.
\end{equation}
Recall further \eqref{varHJ} and \eqref{eq4'} and consider the following approximation for the conditional variance
$$
\sigma^2_{{\bf p}(k)}( n)=\sum_{i\in{\cal I}}\e^{-p_i(k) n}(1-\e^{-p_i (k)  n})-
 ( n)^{-1}\left(\sum_{i\in{\cal I}} n p_i(k) \e^{-p_i (k) n}\right)^2\,.
 $$

\begin{theorem}\label{T2} Notation is the same as in Theorem \ref{T1}. Provided that
$\theta<1$, we have that
$$\lim_{a,b}\frac{
\sigma^2_{{\bf p}(k)}( n)}{\mu_{{\bf p}(k)}(n)}=2^{\theta}-1\qquad \hbox{a.s.}$$
As a consequence, when $(k,n)$ follows the regime \eqref{reg}, 
$$
\frac{N^{(k)}_n-\mu_{{\bf p}(k)}(n)}{\sqrt{(2^{\theta}-1)\mu_{{\bf p}(k)}(n)}}$$
converges in distribution as  to a standard normal variable.
\end{theorem}
It is interesting to note that the parameter $b$ plays no r\^ole in the limits above,
so by a standard  argument based on extraction of sub-sequences, the results still hold under the  weaker requirement that $k,n\to \infty$ such that $k=a\ln n + O(1)$. 
We do not know whether this can be extended to the more general regime when one merely requires that $k\sim a\ln n$. We also underline that the usefulness  of Theorem \ref{T2} is limited by the fact that the centralizing term $\mu_{{\bf p}(k)}(n)$ in the numerator is random,
and only its first order asymptotic \eqref{eq68} is known.

\proof The calculations resemble that in the proof of Theorem \ref{T1}. We start by
 observing that
$$
\sigma^2_{{\bf p}(k)}( n)=\sum_{i\in\N^{k}}Êf(k/a +\ln p_i(k) + c_{k} ) - ( n)^{-1}\left(\sum_{i\in\N^{k}}Êg(k/a +\ln p_i(k) + c_k )\right)^2
$$
with
$$
f(x)= \exp\{-\e^x\}\left(1-\exp\{-\e^x\}\right)\ ,\ 
g(x)= \e^x \exp\{-\e^x\}
$$
and 
$$c_{k}:= \ln n  - k  /a = -b/a + o(1)\,.$$

Easy calculations show that
$$\int_{\R} f(y)\e^{-\theta y}
{\rm d}y= \frac{(2^{\theta}-1)\Gamma(1-\theta)}{\theta} \,,$$
and 
$$\int_{\R} g(y)\e^{-\theta y}
{\rm d}y= \Gamma(1-\theta) \,.$$
Applying Corollary \ref{C1} and  recalling that $\theta=\hbox{\sc m}^{-1}(1/a)$, we deduce that almost surely
\begin{eqnarray*}
\lim_{a,b} \sqrt{k }\e^{-\varphi(\theta)k} \, 
\sigma^2_{{\bf p}(k)}( n)
&=&(2^{\theta}-1)\frac{ \e^{-\theta b/a}\Gamma(1-\theta)}{\theta \sqrt{2\pi \hbox{\sc v}(\theta)}}\,W(\theta)\\
&=& (2^{\theta}-1)\lim_{a,b} \sqrt{k}\,\e^{-\varphi(\theta)k} \mu_{{\bf p}(k)}(n)\,,
\end{eqnarray*}
which is our first  claim. The second then derives 
from Lemma \ref{LHJ}.
\QED

Roughly speaking, the estimation \eqref{eq67} means that for
$a<1/\hbox{\sc m}(1)$, 
$$n\asymp \e^{k/a}\ \Longrightarrow \ N^{(k)}_n\asymp \e^{\varphi(\theta)k}/ \sqrt{k}\,.$$
We shall point out that this is no longer true when 
$a> 1/\hbox{\sc m}(1)$, in other words that a phase transition occurs at the critical value
$a=1/\hbox{\sc m}(1)$. In this direction, we consider the more general  regime for pairs of integers $(k,n)$ such that
\begin{equation}\label{regbis}
k,n\to \infty\quad \hbox{and} \quad 
k\sim a\ln n\,,
\end{equation}
for some fixed $a>0$. Again, when $F(k,n)$ is some function depending on
$k$ and $n$, we shall write
$$\lim_{a}F(k,n)$$
for the limit (whenever it exists)
of $F(k,n)$ when $(k,n)$ follows the regime \eqref{regbis}.

\begin{proposition} \label{P1}  If $1/\hbox{\sc m}(1)< a$, then
$$\lim_{a}\frac{N^{(k)}_{n}}{n}=1\qquad \hbox{a.s.}$$
\end{proposition}

\proof By an argument similar to that in the proof of Theorem \ref{T1}, it suffices to check that 
$$\lim_{a,b} (\alpha n)^{-1} \mu_{{\bf p}(k)}(\alpha n)=1\qquad \hbox{a.s.}$$
provided that $\alpha$ is sufficiently close to $1$. The upperbound
$\mu_{{\bf p}(k)}(\alpha n)\leq \alpha n$ is plain, so we shall focus on the lowerbound.

In this direction, integrating the obvious inequality $(1+\varepsilon)x^{\varepsilon}
\geq 1-\e^{-x}$ for $x\geq 0$ and $0<\varepsilon\leq 1$ arbitrary, we see that
\begin{equation}\label{utile}
1-\e^{-x}\geq x-x^{1+\varepsilon}\,.
\end{equation}
Replacing $x$ by the random variable $\alpha n p_i(k)$ and summing over $i\in\N^k$, we deduce that
$$\mu_{{\bf p}(k)}(\alpha n)\geq \alpha n -(\alpha n)^{1+\varepsilon}\sum_{i\in\N^k}
p_i(k)^{1+\varepsilon}
=\alpha n - (\alpha n)^{1+\varepsilon} \hbox{\sc l}(1+\varepsilon)^kW^{(k)}(1+\varepsilon)\,. $$
As $W^{(k)}(1+\varepsilon)$ is a positive martingale (in the variable $k$), we have $\sup_{k}W^{(k)}(1+\varepsilon)<\infty$ a.s. On the other hand, we have
\begin{eqnarray*}
\ln\left(n^{1+\varepsilon} \hbox{\sc l}(1+\varepsilon)^k\right)
&=& (1+\varepsilon)\ln n + k\ln \hbox{\sc l}(1+\varepsilon)\\
&=& (1+\varepsilon +a \ln \hbox{\sc l}(1+\varepsilon) + o(1))\ln n.
\end{eqnarray*}
Recall that $-\hbox{\sc m}$ is the derivative of $\ln \hbox{\sc l}$ and that
$\ln \hbox{\sc l}(1)=0$ . Since 
$a\hbox{\sc m}(1)>1$, we can chose $\varepsilon>0$ small enough so that 
$a \ln \hbox{\sc l}(1+\varepsilon)<-\varepsilon$.
Then 
$$(\alpha n)^{1+\varepsilon} \hbox{\sc l}(1+\varepsilon)^k=o(n)\,,$$
and we conclude that
$$\liminf_{a,b} (\alpha n)^{-1} \mu_{{\bf p}(k)}(\alpha n)\geq 1\qquad \hbox{a.s.}$$
which completes the proof. \QED

For  $\theta^*\leq 2$, a related argument also enables us to estimate the lowest generation 
at which all $n$ balls fall into different boxes. More generally, recall that for every integer $j$, $\bar N^{(k)}_{n,j}$ denotes the number of boxes at generation $k$ which contain more than $j$ balls when $n$ balls have been thrown. 
Note that this quantity increases with the number of balls $n$ and decreases with the generation $k$. Define
$$\zeta_{n,j}:=\min\{k\in\N: \bar N^{(k)}_{n,j}=0\}\,.$$

\begin{proposition} \label{P2} 
If $\theta^*\leq j+1$,
then we have 
$$\lim_{n\to\infty} \frac{\zeta_{n,j}}{\ln n}= 1/\hbox{\sc m}_*\qquad \hbox{a.s.}$$
\end{proposition}

\proof 
We first consider the randomized version with a total number of balls ${\tt n}_x$ 
which has the Poisson law with parameter $x$. We write $\bar {\tt N}_{x,j}^{(k)}$ for the number of  boxes at generation $k$ occupied by more than $j$ balls, and consider its  conditional expectation  given the random probability measure ${\bf p}(k)$, which has been computed in \eqref{eq45}. Using  the assumption that $\theta^*\leq j+1$
at the second line below, we get 
\begin{eqnarray*}
\E(\bar {\tt N}_{x,j}^{(k)}\mid {\bf p}(k)) = \bar 
\mu_{{\bf p}(k),j}(x)&=&\sum_{i\in\N^k} \left(  1-\e^{-p_i (k)x}
\left( 1+\cdots +\frac{(x p_i(k))^j}{j!}\right)\right)\\
&\leq & c(\theta^*) \sum_{i\in\N^k} \left(  x p_i(k)\right)^{\theta^*}\,,
\end{eqnarray*}
where 
$$c(\theta^*)=\max_{y>0} y^{-\theta^*}
\left(  1-\e^{-y}\left( 1+\cdots +\frac{y^j}{j!}\right)\right)$$ is some finite constant.
Observe that the preceding upperbound can be expressed as
$$ c(\theta^*) x^{\theta^*} 
\hbox{\sc l}(\theta^*)^kW^{(k)}(\theta^*)\,,$$
and recall that  $W^{(k)}(\theta^*)$ is a martingale.
The unconditional expectation
$\E(\bar{\tt N}_{x,j}^{(k)})$ can thus be bounded from above by $c(\theta^*)x^{\theta^*} 
\hbox{\sc l}(\theta^*)^k$.

We next pick $a>1/\hbox{\sc m}_*$ and take  $x=\e^{k/a}$.
Recall  that 
$$-\frac{\ln \hbox{\sc l}(\theta^*)}{\theta^*}=\hbox{\sc m}(\theta^*)
=\hbox{\sc m}_*\,.$$
We thus have
$$\ln\left( x^{\theta^*} 
\hbox{\sc l}(\theta^*)^k\right)= k\theta^* \left(\frac{1}{a}-\hbox{\sc m}_*\right)\,,$$
and as a consequence
\begin{equation}\label{moment}
\sum_{k\in\N}\E(\bar{\tt N}_{\e^{k/a},j}^{(k)})<\infty.
\end{equation}

Then chose any $a'>a$ and recall that the Poisson process fulfills
$$\lim_{k_0\to \infty}\P\left({\tt n}_{\e^{k/a}}>\lfloor \e^{(k+1)/a'} \rfloor \hbox{ for all }k\geq k_0\right)
=1\,.$$
We deduce from \eqref{moment} and an argument of monotonicity that 
$$\P\left(\bar{N}_{\lfloor \e^{(k+1)/a'} \rfloor,j}^{(k)}=0 \hbox{ for all integers $k$ sufficiently large} \right)=1.$$
Observing that for every integer $n$ and $k=\lfloor a' \ln n\rfloor$:
$$\zeta_{n,j}> k\ \Longrightarrow\ \bar N^{(k)}_{n,j}\geq 1
\Longrightarrow\ \bar N^{(k)}_{\lfloor \e^{(k+1)/a'} \rfloor,j}\geq 1\,,$$
and therefore
$$\limsup_{n\to\infty} \frac{\zeta_{n,j}}{\ln n}\leq  a'\qquad \hbox{a.s.}$$
As $a'$ can be chosen arbitrarily close to $1/\hbox{\sc m}_*$, we conclude that
$$\limsup_{n\to\infty} \frac{\zeta_{n,j}}{\ln n}\leq  1/\hbox{\sc m}_*\qquad \hbox{a.s.}$$
Finally, the converse bound 
$$\liminf_{n\to\infty} \frac{\zeta_{n,j}}{\ln n}\geq  1/\hbox{\sc m}_*\qquad \hbox{a.s.}$$ derives easily  from Theorem \ref{T1}. \QED

\subsection{Interpretation in terms of homogeneous fragmentations}

The initial motivation for this work was to gain insight on certain asymptotic regimes for
 homogeneous fragmentations. Roughly, the latter form coherent families of  natural Markov processes with values in the space of partitions of finite sets, such that these random partitions are refining as time passes. They are closely related to multiplicative cascades of random probability measures and to the occupancy scheme;  we start by giving precise definitions.
 
A non-empty set $B$ of positive integers  is called a block, and
a partition of $B$  is a denumerable family $\pmb{\pi}=\{\pi_1,\pi_2,\ldots\}$ of pairwise disjoint 
sub-blocks of $B$ such that $\cup\pi_i=B$. We write $\p_B$ for the set of partitions of $B$ and endow  $\p_B$ with a natural partial order~: one says that a partition $\pmb{\pi}$ is finer than another partition  $\pmb{\pi}'$ and then write $\pmb{\pi}\preceq \pmb{\pi}'$ if and only if each block $\pi_i$ of $\pmb{\pi}$ is contained into some block $\pi'_j$ of $\pmb{\pi}'$. 
A sequence $(\pmb{\pi}(k), k\geq 0)$ of partitions  is called {\it nested} (or, sometimes also, refining) if $\pmb{\pi}(k+1)$ is finer than $\pmb{\pi}(k)$ for every integer $k\geq 0$.  

The occupancy scheme produces naturally random partitions. 
Typically, we consider some discrete probability measure ${\bf p}\in{\rm Prob}_{\cal I}$
and the corresponding family of boxes, and we label balls by integers. Every block $B$ can then be splitted into sub-blocks that correspond to the labels of the  balls which occupy the same box. This provides a random partition of $B$, say $\pmb{\pi}^{\bf p}_B$.
The latter is  {\it exchangeable}, in the sense that its distribution is invariant under the natural action of permutations of $B$.  Note that when $B$ is infinite, in particular when $B=\N$,  
$\pmb{\pi}^{\bf p}_B$ has no singletons a.s. A fundamental theorem due to Kingman \cite{King} (see Theorem 2.1 in \cite{RFCP}) claims that any exchangeable random partition of $\N$ which has no singletons a.s. has the same distribution as the partition that results from some randomized version of the occupancy scheme, i.e. for which ${\bf p}$ is now a random discrete probability measure.
The assumption of absence of singletons can be dropped provided that one allows ${\bf p}$ to be defective, i.e. to be only a sub-probability measure.

We now  consider again a sequence $({\bf p}(k), k\in\N)$ of random discrete probability measures which is associated to some multiplicative cascade as in \eqref{eq3}, 
and for every  $k\in\N$, we denote by $\pmb{\Pi}(k)$ the random partition of $\N$ 
induced as above by the occupancy scheme at generation $k$.
Then $\pmb{\Pi}=(\pmb{\Pi}(k), k\geq0)$ is a  nested sequence of exchangeable random partitions of $\N$,  which is Markovian. We call $\pmb{\Pi}$ a {\it homogeneous fragmentation chain}. Its transition probabilities inherit the branching property from the multiplicative structure of the cascade; see Propositions 1.2 and 1.3 in \cite{RFCP}.

For any of blocks $B$ and $B'$ with $B'\subseteq B$, the restriction to $B'$ yields a natural projection
$\pmb{\pi} \to \pmb{\pi}_{\mid B'}$ from $\p_B$ to $\p_{B'}$. The partial order $\preceq $ is clearly compatible with restrictions, in the sense that $\pmb{\pi}\preceq \pmb{\pi}' \Rightarrow \pmb{\pi}_{\mid B'}\preceq \pmb{\pi}'_{\mid B'}$. Thus, if $(\pmb{\pi}(k), k\geq 0)$ is a nested sequence of partitions of a block $B$ and if $B'\subseteq B$ is a smaller block,
then the sequence $(\pmb{\pi}_{\mid B'}(k), k\geq 0)$ of partitions restricted to  $B'$ is again nested. A simple but nonetheless important fact is that the Markov property of a homogeneous fragmentation chain $\pmb{\Pi}$ is preserved by restriction, in the sense that for every block $B\subseteq \N$, the nested sequence of  partitions of $B$, $\pmb{\Pi}_{\mid B}=(\pmb{\Pi}_{\mid B}(k), k\geq0)$, is still Markovian;
see Lemma 3.4 in \cite{RFCP} for a sharper statement.

Recapitulating, the occupancy scheme enables us to associate to any random multiplicative cascade of discrete probability measures a homogeneous fragmentation chain  $\pmb{\Pi}$. In turn, the latter
provided a nested sequence of random partitions of an arbitrary block $B\subseteq \N$, 
$\pmb{\Pi}_{\mid B}=(\pmb{\Pi}_{\mid B}(k), k\geq0)$, which is Markovian. Further, these Markov chains are coherent, in the sense that if $B'\subseteq B$, then $\pmb{\Pi}_{\mid B'}
= (\pmb{\Pi}_{\mid B})_{\mid B'}$. Roughly speaking,
the statements in the preceding Section provide information about the asymptotic regimes for a homogeneous fragmentation of a finite set, when both the size $n$ of that set and the time $k$ at which the fragmentation process is observed 
tend to infinity. For example, Theorem \ref{T1} is a limit theorem 
 for the number of components with a fixed size (like singletons, pairs, etc.),
whereas Proposition \ref{P2} specifies the asymptotic behavior of the shattering time, that is the first instant at which the fragmentation process of a finite block reaches its absorbing state, i.e. the partition of that block into singletons. 
Finally, we also mention that, for the sake of simplicity, we have only discussed here fragmentation processes in discrete time; however our results can be shifted to homogeneous fragmentations in continuous time, using discretization techniques  similar to those in developed in \cite{BeRou}. 

\end{section}

\vskip 4mm
\noindent{\large \bf Appendix : Proof of the large deviations behavior}

We finally proceed to the proof of Corollary \ref{C1}.
 When $f$ is continuous with compact support, the claim follows from Lemma \ref{L2} by approximating $f$ with 
step functions. See e.g. Corollary 4 of \cite{Big2} and Theorem 3 of Stone \cite{St} for slightly stronger statements 
in terms of directly Riemann integrable functions with compact support.
All that is needed to extend this to continuous functions with unbounded support  (which fulfill the conditions of the statement) is to establish the following : If we define  for some  fixed $\alpha>0$ and $\beta >\theta$
$$g_+(x)=  {\bf 1}_{\{x>0\}}x^{\alpha }\,,\quad
\hbox{and}Ê\quad g_-(x)={\bf 1}_{\{x<0\}} \e^{\beta x}\,,$$
then,
\begin{equation}\label{eq22}
\sup_{k\in\N}  \sqrt{k }\e^{-\varphi(\theta)k}\sum_{i\in\N^k}Êg_{\pm}(k \hbox{\sc m}(\theta) +\ln p_i(k))<\infty \qquad\hbox{a.s.}
\end{equation}
Indeed, for every function $f$ that fulfills the hypotheses of the statement and for every integer $\ell\geq 1$, we can find a continuous function $f_{\ell}$ with compact support such that $f_{\ell}\leq f \leq f_{\ell} + \ell^{-1}(g_++g_-)$, and then \eqref{eq22} enables us to conclude the proof by  a standard argument.

\noindent{\bf Proof of \eqref{eq22} for $g_+$:}\hskip10pt
We write
$$
\sum_{i\in\N^k}Êg_{+}(k \hbox{\sc m}(\theta) +\ln p_i(k))
= \int_{\R_+}Êg_+(k \hbox{\sc m}(\theta) -y) Z^{(k)}(  {\rm d}y)
= \int_{[0,k \hbox{\sc m}(\theta)]}Ê(k \hbox{\sc m}(\theta) -y)^{\alpha} Z^{(k)}(  {\rm d}y)\,.$$
Then we pick $\theta'\in]\theta,\theta^*[$ sufficiently close to $\theta$ (as this will be explained in the sequel) and split the last integral at $k \hbox{\sc m}(\theta')$ to get
$$  \int_{[0,k \hbox{\sc m}(\theta')]}Ê(k \hbox{\sc m}(\theta) -y)^{\alpha} Z^{(k)}(  {\rm d}y)
+\int_{]k \hbox{\sc m}(\theta'),k \hbox{\sc m}(\theta)]}Ê(k \hbox{\sc m}(\theta) -y)^{\alpha} Z^{(k)}(  {\rm d}y)\,.
$$
For the first integral, there is the obvious bound
$$\int_{[0,k \hbox{\sc m}(\theta')]}Ê(k \hbox{\sc m}(\theta) -y)^{\alpha} Z^{(k)}(  {\rm d}y) \leq
(k \hbox{\sc m}(\theta))^{\alpha} Z^{(k)}([0,k \hbox{\sc m}(\theta')])\,.$$
Observe from Markov inequality that  $Z^{(k)}([0,a])
\leq \e^{\theta a} \hbox{\sc l}(\theta)^kW^{(k)}(\theta)$ for every $a>0$,
so the preceding quantity can be bounded from above by 
$$ (k \hbox{\sc m}(\theta))^{\alpha}\e^{  \theta  k  \hbox{\sc m}(\theta')} \hbox{\sc l}(\theta)^kW^{(k)}(\theta).$$
Recall that $\varphi(\theta) = \ln\hbox{\sc l}(\theta) + \theta \hbox{\sc m}(\theta)$
and that $\hbox{\sc m}(\theta')<\hbox{\sc m}(\theta)$.
As a consequence
$$(k \hbox{\sc m}(\theta))^{\alpha}\e^{ \theta k   \hbox{\sc m}(\theta')}
 \hbox{\sc l}(\theta)^k
= o(\e^{k\varphi(\theta)})/\sqrt k \,,\qquad k\to\infty\,,$$
and since the martingale $W^{(k)}(\theta)$ remains bounded a.s., we conclude that
\begin{equation}\label{eq24}
\lim_{k\to\infty}\sqrt k \e^{-k\varphi(\theta)} \int_{[0,k \hbox{\sc m}(\theta')]}Ê(k \hbox{\sc m}(\theta) -y)^{\alpha} Z^{(k)}(  {\rm d}y) = 0\qquad \hbox{a.s.}
\end{equation}

For the second integral, we start from the bound
\begin{eqnarray*}
& &\int_{]k \hbox{\sc m}(\theta'),k \hbox{\sc m}(\theta)]}Ê(k \hbox{\sc m}(\theta) -y)^{\alpha} Z^{(k)}(  {\rm d}y)\\
&\leq&
\sum_{0\leq j \leq k( \hbox{\sc m}(\theta)- \hbox{\sc m}(\theta'))}j^{\alpha}Z^{(k)}([ k \hbox{\sc m}(\theta) -j,  k \hbox{\sc m}(\theta)-j +1]).
\end{eqnarray*}
For indices $0\leq j \leq k( \hbox{\sc m}(\theta)- \hbox{\sc m}(\theta'))
$, we define 
\begin{equation}\label{eq23}
\theta_{k,j}= \hbox{\sc m}^{-1}(\hbox{\sc m}(\theta)-j/k)  \in [\theta,\theta'] \,,
\end{equation}
so that
$$ k \hbox{\sc m}(\theta_{k,j})= k \hbox{\sc m}(\theta) -j\,.$$
We observe that
\begin{eqnarray*} 
Z^{(k)}([ k \hbox{\sc m}(\theta) -j,  k \hbox{\sc m}(\theta)-j +1])
&=&  \hbox{\sc l}(\theta_{k,j})^k\int_{k \hbox{\sc m}(\theta)-j}^{k \hbox{\sc m}(\theta) -j +1} 
 \e^{\theta_{k,j}y}Z^{(k)}_{\theta_{k,j}}(  {\rm d}y)
\\
 &\leq& \e^{
\theta_{k,j}( k \hbox{\sc m}(\theta)-j+1)} \hbox{\sc l}(\theta_{k,j})^k
Z^{(k)}_{\theta_{k,j}}([ k \hbox{\sc m}(\theta) -j,  k \hbox{\sc m}(\theta) -j +1])\\
&\leq & \e^{(1-j)\theta}
\left(\e^{\theta_{k,j}  \hbox{\sc m}(\theta)} \hbox{\sc l}(\theta_{k,j})\right)^k
Z^{(k)}_{\theta_{k,j}}([ k \hbox{\sc m}(\theta_{k,j}),  k \hbox{\sc m}(\theta_{k,j})+1])\,,
\end{eqnarray*}
where for the last inequality, we use the fact that 
$\theta_{k,j}\geq \theta$.

Recall  that $\ln \hbox{\sc l}$ is convex and has derivative $-\hbox{\sc m}$.
Since $\theta\leq \theta_{k,j}\leq \theta'$, we have
$$\ln \hbox{\sc l}(\theta_{k,j})\leq  \ln \hbox{\sc l}(\theta) -
\hbox{\sc m}(\theta')(\theta_{k,j}-\theta)\,,$$
and since $\varphi(\theta) = \ln\hbox{\sc l}(\theta) + \theta \hbox{\sc m}(\theta)$, this yields
$${\theta_{k,j}  \hbox{\sc m}(\theta)} + \ln\hbox{\sc l}(\theta_{k,j})
\leq \varphi(\theta) + (\theta_{k,j}-\theta)(\hbox{\sc m}(\theta)-\hbox{\sc m}(\theta')).$$
As the mean function $\hbox{\sc m}$ is locally a regular diffeomorphism, we see from \eqref{eq23} that there is some finite
constant $C$ (which is independent of $k$ and $j$) such that $\theta_{k,j}-\theta\leq C j/k$. Further, we also have that 
$\hbox{\sc m}(\theta)-\hbox{\sc m}(\theta')\leq \theta/2C$ provided that we choose $\theta'$ sufficiently close to $\theta$.
Then 
$$\left(\e^{\theta_{k,j}  \hbox{\sc m}(\theta)} \hbox{\sc l}(\theta_{k,j})\right)^k \leq \e^{j\theta/2} \exp(k\varphi(\theta))\,,$$
and thus 
$$\sum_{0\leq j \leq k( \hbox{\sc m}(\theta)- \hbox{\sc m}(\theta'))} j^{\alpha}\e^{(1-j)\theta}
\left(\e^{\theta_{k,j}  \hbox{\sc m}(\theta)} 
\hbox{\sc l}(\theta_{k,j})\right)^k 
= O\left(\exp(k\varphi(\theta))\right)\,.$$

On the other hand, we know from Lemma \ref{L2} that there exists an a.s. finite random variable $\xi$ such that
$$Z^{(k)}_{\theta_{k,j}}([ k \hbox{\sc m}(\theta_{k,j}),  k \hbox{\sc m}(\theta_{k,j})+1])\leq k^{-1/2}\xi\,.$$
Putting the pieces together, we conclude that
$$\sup_{k\in\N} \sqrt k   \e^{-k\varphi(\theta)} 
\int_{]k \hbox{\sc m}(\theta'),k \hbox{\sc m}(\theta)]}Ê(k \hbox{\sc m}(\theta) -y)^{\alpha} Z^{(k)}(  {\rm d}y) <\infty\qquad 
\hbox{a.s.}
$$
Combining with \eqref{eq24}, we have thus checked that \eqref{eq22} does hold for $g_+$. \QED

The proof of  the bound \eqref{eq22} for $g_-$ follows a similar route; however it may be useful to spell out the main steps.

\noindent{\bf Proof of \eqref{eq22} for $g_-$:}\hskip10pt
We start with 
$$
\sum_{i\in\N^k}Êg_{-}(k \hbox{\sc m}(\theta) +\ln p_i(k))
= \int_{[0,\infty[}Ê\e^{-\beta y} Z^{(k)}( k \hbox{\sc m}(\theta)+ {\rm d}y).
$$
We pick $\theta'\in]\theta_*,\theta[$ sufficiently close to $\theta$ and split this integral at $ k( \hbox{\sc m}(\theta')- \hbox{\sc m}(\theta))$.

For indices $0\leq j < k( \hbox{\sc m}(\theta')- \hbox{\sc m}(\theta))$, we introduce
\begin{equation}\label{eq25}
\theta_{k,j}= \hbox{\sc m}^{-1}(\hbox{\sc m}(\theta)+j/k)  \in [\theta',\theta] \,,
\end{equation}
so that
$$ k \hbox{\sc m}(\theta_{k,j})= k \hbox{\sc m}(\theta) + j\,.$$
We observe that
\begin{eqnarray*} & &
\int_{[j,j+1[}Ê\e^{-\beta y} Z^{(k)}( k \hbox{\sc m}(\theta)+ {\rm d}y)\\
&=& \hbox{\sc l}(\theta_{k,j})^k
\int_{[j,j+1[}Ê\e^{-\beta y} \e^{ \theta_{k,j} (k \hbox{\sc m}(\theta)+ y)}
Z^{(k)}_{\theta_{k,j}}( k \hbox{\sc m}(\theta)+ {\rm d}y)\\
 &\leq&  \hbox{\sc l}(\theta_{k,j})^k \e^{-j(\beta -
\theta_{k,j})}\e^{ k \theta_{k,j}\hbox{\sc m}(\theta)}
Z^{(k)}_{\theta_{k,j}}([ k \hbox{\sc m}(\theta) +j,  k \hbox{\sc m}(\theta) +j+1])\\
&\leq & \e^{-j(\beta -\theta)}
\left(\e^{\theta_{k,j}  \hbox{\sc m}(\theta)} \hbox{\sc l}(\theta_{k,j})\right)^k
Z^{(k)}_{\theta_{k,j}}([ k \hbox{\sc m}(\theta_{k,j}),  k \hbox{\sc m}(\theta_{k,j})+1])\,,
\end{eqnarray*}
where for the last inequality, we use the fact that 
$\theta_{k,j}\leq \theta$.

Because $\ln \hbox{\sc l}$ is convex, has derivative $-\hbox{\sc m}$ and
 $\theta'\leq \theta_{k,j}\leq \theta$, there is the inequality
$$\ln \hbox{\sc l}(\theta_{k,j})\leq  \ln \hbox{\sc l}(\theta) +
\hbox{\sc m}(\theta')(\theta-\theta_{k,j})\,,$$
which yields
$${\theta_{k,j}  \hbox{\sc m}(\theta)} + \ln\hbox{\sc l}(\theta_{k,j})
\leq \varphi(\theta) + (\theta-\theta_{k,j})(\hbox{\sc m}(\theta')-\hbox{\sc m}(\theta)).$$
We see from \eqref{eq25} that there is some finite
constant $C$ (which is independent of $k$ and $j$) such that $\theta-\theta_{k,j}\leq C j/k$, and that 
$\hbox{\sc m}(\theta')-\hbox{\sc m}(\theta)\leq (\beta-\theta)/2C$ provided that we choose $\theta'$ sufficiently close to $\theta$.
Then 
$$\left(\e^{\theta_{k,j}  \hbox{\sc m}(\theta)} \hbox{\sc l}(\theta_{k,j})\right)^k \leq \e^{j(\beta-\theta)/2} \exp(k\varphi(\theta))\,,$$
and thus 
$$\sum_{0\leq j < k( \hbox{\sc m}(\theta')- \hbox{\sc m}(\theta))} \e^{-j(\beta-\theta)}
\left(\e^{\theta_{k,j}  \hbox{\sc m}(\theta)} \hbox{\sc l}(\theta_{k,j})\right)^k = O\left(\exp(k\varphi(\theta))\right)\,.$$
Further Lemma \ref{L2} ensures the existence of an a.s. finite random variable $\xi$ such that
$$Z^{(k)}_{\theta_{k,j}}([ k \hbox{\sc m}(\theta_{k,j}),  k \hbox{\sc m}(\theta_{k,j})+1])\leq k^{-1/2}\xi\,,$$
and then we can conclude that
\begin{equation}\label{eq26}
\sup_{k\in\N} \sqrt k   \e^{-k\varphi(\theta)} 
\int_{[0, k( \hbox{\sc m}(\theta')- \hbox{\sc m}(\theta))[}
Ê\e^{-\beta y} Z^{(k)}( k \hbox{\sc m}(\theta)+ {\rm d}y)
<\infty\qquad 
\hbox{a.s.}
\end{equation}

For the remaining integral, we  
note that
\begin{eqnarray*}
& &\int_{[k( \hbox{\sc m}(\theta')- \hbox{\sc m}(\theta)),\infty[}
Ê\e^{-\beta y} Z^{(k)}( k \hbox{\sc m}(\theta)+ {\rm d}y) \\
&\leq&
\e^{-(\beta-\theta) k( \hbox{\sc m}(\theta')- \hbox{\sc m}(\theta))}
\int_{[k( \hbox{\sc m}(\theta')- \hbox{\sc m}(\theta)),\infty[}
Ê\e^{-\theta y} Z^{(k)}( k \hbox{\sc m}(\theta)+ {\rm d}y) \\
&\leq&
\e^{-(\beta-\theta) k( \hbox{\sc m}(\theta')- \hbox{\sc m}(\theta))}
\e^{k\theta \hbox{\sc m}(\theta)} \hbox{\sc l}(\theta)^k W^{(k)}(\theta)\,.
 \end{eqnarray*}
 We readily deduce that
 $$\lim_{k\to\infty} \sqrt k   \e^{-k\varphi(\theta)} 
\int_{[k( \hbox{\sc m}(\theta')- \hbox{\sc m}(\theta)),\infty[}
Ê\e^{-\beta y} Z^{(k)}( k \hbox{\sc m}(\theta)+ {\rm d}y) = 0\qquad 
\hbox{a.s.,}$$
and combining with \eqref{eq26}, this establishes that \eqref{eq22} holds for $g_-$. \QED

\vskip 1cm
\noindent {\bf Acknowledgment.} I would like to thank an anonymous referee for 
having carefully checked of the first draft of this work, and especially for pointing at some errors.


\begin{thebibliography}{99}

\bibitem{Bahadur}
R. R. Bahadur (1960). On the number of distinct values in a large sample
from an infinite discrete distribution. {\sl Proc. Nat. Inst. Sci. India Part A} {\bf 26}, 67-75.

\bibitem{RFCP} J. Bertoin  (2006). {\em Random Fragmentation and Coagulation Processes}. Cambridge University Press, Cambridge.
 
 \bibitem{BeRou} J. Bertoin and  A. Rouault (2005). Discretization methods for homogeneous fragmentations. 
{\sl  J. London Math. Soc.} {\bf 72},  91-109. 

\bibitem{Big1}  J. D. Biggins (1977). Martingale convergence in the branching
random
walk. {\sl  J. Appl. Probability } {\bf 14}, 25-37.


\bibitem{Big2} J. D.  Biggins  (1992). Uniform convergence of martingales in the
branching random walk.  {\sl  Ann. Probab.} {\bf 20}, 137-151.

\bibitem{Bou} S. Boucheron and D. Gardy (1997). An urn model from learning theory. 
{\sl Random Struct. Algorithms \bf 10}, 43-69.

\bibitem{Bun} J. Bunge and M. Fitzpatrick (1993). Estimating the number of species : a review. {\sl J. Am. Statist. Assoc. \bf 88}, 364-373.

\bibitem{Dutko} 
M. Dutko (1989). Central limit theorems for infinite urn models. {\sl Ann.
Probab. \bf 17}, 1255-1263.

\bibitem{Gardy}  D. Gardy (2002). Occupancy urn models in the analysis of algorithms. 
{\it J. Stat. Plann. Inference \bf 101}, 95-105.

\bibitem{GHP}
A. Gnedin, B. Hansen and J. Pitman (2007). Notes on the occupancy problem
with infinitely many boxes: general asymptotics and power laws. {\sl Probability Surveys \bf 4},
146-171. 

\bibitem{GPY}
A. Gnedin, J. Pitman and M. Yor (2006). Asymptotic laws for compositions
derived from transformed subordinators. {\sl Ann. Probab. \bf 34}, 468-492.

\bibitem{HeG} F. He and  K. J. Gaston (2000). 
Estimating species abundance from occurrence. {\sl 
The American Naturalist \bf 156-5}, 553-559.

\bibitem{Holst} L. Holst (1986). On birthday, collectors', occupancy and other classical urn problems. {\sl Int. Stat. Rev. \bf 54}, 15-27.

\bibitem{HJ} H.-K. Hwang and S. Janson (2007+).
Local limit theorems for finite and infinite urn models. To appear in {\sl Ann. Probab.}

\bibitem{JK}
N. L. Johnson and S. Kotz (1977). {\sl Urn Models and Their Application. An
Approach to Modern Discrete Probability Theory}. John Wiley \& Sons, New
York-London-Sydney.

\bibitem{Karlin}
S. Karlin (1967). Central limit theorems for certain infinite urn schemes. {\sl J.
Math. Mech. \bf 17}, 373-401.

\bibitem{King} J. F. C. Kingman(1982).
 The coalescent.
{\sl  Stochastic Process. Appl.} {\bf 13}, 235-248.

\bibitem{Kolchin}
V. F. Kolchin, B. A. Sevast'yanov and V. P. Chistyakov (1978). {\sl Random
Allocations.} John Wiley \& Sons, New
York-London-Sydney.

\bibitem{Liu} Q. S. Liu (2000). On generalized multiplicative cascades.
{\sl Stochastic Process. Appl. \bf 86}, 263-286.

\bibitem{St} C. Stone (1967). On local and ratio limit theorems. {\sl
Proc. 5th Berkeley Sympos. math. Statist. Probab. \bf 2},  217-224.
\end{thebibliography}
\end{document}